\RequirePackage{fix-cm}
\documentclass[twocolumn]{svjour3}          
\smartqed  
\usepackage{graphicx}
\usepackage{epstopdf}
\usepackage{amssymb}
\usepackage{amsmath}
\usepackage{enumerate}
\usepackage[none]{hyphenat}
 \usepackage{mathptmx}      
\spnewtheorem{algorithm}{Algorithm}{\bf}{\rm}
\begin{document}

\title{On linear and quadratic Lipschitz bounds for twice continuously differentiable functions
}


\author{Gene A. Bunin, Gr\'egory Fran\c cois, Dominique Bonvin 
}


\institute{Laboratoire d'Automatique, Ecole Polytechnique F\' ed\' erale de Lausanne, Lausanne, CH-1015 \\
              \email{gene.a.bunin@ccapprox.info}           
}

\date{Submitted: \today}

\maketitle

\begin{abstract}
Lower and upper bounds for a given function are important in many mathematical and engineering contexts, where they often serve as a base for both analysis and application. In this short paper, we derive piecewise linear and quadratic bounds that are stated in terms of the Lipschitz constants of the function and the Lipschitz constants of its partial derivatives, and serve to bound the function's evolution over a compact set. While the results follow from basic mathematical principles and are certainly not new, we present them as they are, from our experience, very difficult to find explicitly either in the literature or in most analysis textbooks.
\end{abstract}

\section{Overview}
\label{summary}

We consider a function $f : \mathbb{R}^n \rightarrow \mathbb{R}$ of the variables $x \in \mathbb{R}^n$ that is twice continuously differentiable ($C^2$) over an open set containing the compact set $\mathcal{X}$. Because $f$ is $C^2$ over $\mathcal{X}$, its first and second derivatives on this set exist and must be bounded by the Lipschitz constants

\vspace{-4mm}
\begin{equation}\label{eq:lip1}
\underline \kappa_i < \frac{\partial f}{\partial x_i} \Big |_{x} < \overline \kappa_i, \hspace{15mm} i = 1,...,n
\end{equation}

\vspace{-4mm}
\begin{equation}\label{eq:lip2}
\underline M_{ij} < \frac{\partial^2 f}{\partial x_j \partial x_i} \Big |_{x} < \overline M_{ij}, \hspace{6mm}  i,j = 1,...,n
\end{equation}

\noindent for all $x \in \mathcal{X}$.

The evolution of $f$ between any two points $x_a, x_b \in \mathcal{X}$ may then be bounded as

\vspace{-4mm}
\begin{equation}\label{eq:bound1L}
f(x_b) - f(x_a) \geq \displaystyle \sum_{i=1}^n \mathop {\min} \left[ \begin{array}{l} \underline \kappa_i (x_{b,i} - x_{a,i}), \\ \overline \kappa_i (x_{b,i} - x_{a,i}) \end{array} \right],
\end{equation}

\vspace{-4mm}
\begin{equation}\label{eq:bound1U}
f(x_b) - f(x_a) \leq \displaystyle \sum_{i=1}^n \mathop {\max} \left[ \begin{array}{l} \underline \kappa_i (x_{b,i} - x_{a,i}), \\ \overline \kappa_i (x_{b,i} - x_{a,i}) \end{array} \right],
\end{equation}

\noindent where $x_{a,i}$ and $x_{b,i}$ denote the $i^{\rm th}$ elements of the vectors $x_a$ and $x_b$, respectively.

The bounds (\ref{eq:bound1L}) and (\ref{eq:bound1U}) are piecewise linear in $x$. Alternatively, one may also use the piecewise quadratic bounds

\vspace{-4mm}
\begin{equation}\label{eq:bound2L}
\begin{array}{l}
f(x_b) - f(x_a) \geq \nabla f (x_a)^T (x_b - x_a) + \\
\hspace{10mm}\displaystyle \frac{1}{2} \sum_{i=1}^n \sum_{j=1}^n \mathop {\min} \left[ \begin{array}{l} \underline M_{ij} (x_{b,i} - x_{a,i})(x_{b,j} - x_{a,j}), \\ \overline M_{ij} (x_{b,i} - x_{a,i})(x_{b,j} - x_{a,j}) \end{array} \right],
\end{array}
\end{equation}

\vspace{-4mm}
\begin{equation}\label{eq:bound2U}
\begin{array}{l}
f(x_b) - f(x_a) \leq \nabla f (x_a)^T (x_b - x_a) + \\
\hspace{10mm}\displaystyle \frac{1}{2} \sum_{i=1}^n \sum_{j=1}^n \mathop {\max} \left[ \begin{array}{l} \underline M_{ij} (x_{b,i} - x_{a,i})(x_{b,j} - x_{a,j}), \\ \overline M_{ij} (x_{b,i} - x_{a,i})(x_{b,j} - x_{a,j}) \end{array} \right],
\end{array}
\end{equation}

\noindent which are locally less conservative but also require more knowledge in the form of both the gradient and the Lipschitz constants on the partial derivatives of $f$. While one may generalize this pattern to even higher orders, we will content ourselves with the linear and quadratic cases as we believe these to be sufficient for most applications -- see, however, \cite{Cartis:13} for a discussion of the cubic case.

\section{Derivation of the Linear Bounds}

To limit our analysis to a single dimension, we will consider the line segment between $x_a$ and $x_b$. The following one-dimensional parameterization is used:

\vspace{-4mm}
\begin{equation}\label{eq:Qproof1}
\hat f (\gamma) = f(x (\gamma)),
\end{equation}

\noindent with $x (\gamma) = x_a + \gamma (x_b - x_a), \; \gamma \in [0,1]$. As $f$ is $C^2$, it follows that $\hat f$ is as well, which allows us to use the Taylor series expansion between $\gamma = 0$ and $\gamma = 1$, together with the mean-value theorem \cite{Korn:00}, to state:

\vspace{-4mm}
\begin{equation}\label{eq:Lproof2}
\begin{array}{l}
\hat f (1) = \hat f (0) + \displaystyle \frac{d \hat f}{d \gamma} \Big |_{\tilde \gamma} 
\end{array},
\end{equation}

\noindent for some $\tilde \gamma \in (0,1)$. We proceed to define the first-order derivative in terms of the original function $f$. To do this we apply the chain rule:

\vspace{-4mm}
\begin{equation}\label{eq:Lproof3}
\displaystyle \frac{d \hat f}{d \gamma} \Big |_{\gamma} = \displaystyle \mathop {\sum} \limits_{i = 1}^{n} \frac{\partial f}{\partial x_i} \Big |_{x (\gamma)} \frac{d x_i}{d \gamma} \Big |_{\gamma} = \displaystyle \mathop {\sum} \limits_{i = 1}^{n} \frac{\partial f}{\partial x_i} \Big |_{x (\gamma)} (x_{b,i}-x_{a,i}).
\end{equation}

Noting that $\hat f (0) = f(x_a)$ and $\hat f (1) = f(x_b)$, one may substitute (\ref{eq:Lproof3}) into (\ref{eq:Lproof2}) to obtain

\vspace{-4mm}
\begin{equation}\label{eq:Lproof4}
\begin{array}{l}
f (x_b) = f (x_a) + \displaystyle \mathop {\sum} \limits_{i = 1}^{n} \frac{\partial f}{\partial x_i} \Big |_{x (\tilde \gamma)} (x_{b,i}-x_{a,i}) 
\end{array}.
\end{equation}

Because $ x(\tilde \gamma) \in \mathcal{X}$, we may use (\ref{eq:lip1}) to bound the individual summation components as

\vspace{-4mm}
\begin{equation}\label{eq:mvalbound1}
\begin{array}{l}
x_{b,i} - x_{a,i} \geq 0 \Leftrightarrow \\
\hspace{2mm} \underline \kappa_i (x_{b,i} - x_{a,i}) \leq \displaystyle \frac{\partial f}{\partial x_i} \Big |_{x(\tilde \gamma)} (x_{b,i} - x_{a,i} ) \leq \overline \kappa_i (x_{b,i} - x_{a,i}), \\
x_{b,i} - x_{a,i} \leq 0 \Leftrightarrow \\
\hspace{2mm} \overline \kappa_i (x_{b,i} - x_{a,i}) \leq \displaystyle \frac{\partial f}{\partial x_i} \Big |_{x(\tilde \gamma)} (x_{b,i} - x_{a,i} ) \leq \underline \kappa_i (x_{b,i} - x_{a,i}),
\end{array}
\end{equation}

\noindent or, to account for both cases, as

\vspace{-4mm}
\begin{equation}\label{eq:mvalbound1pw}
\begin{array}{l}
\mathop {\min} \left[ \begin{array}{l} \underline \kappa_i (x_{b,i} - x_{a,i}), \\ \overline \kappa_i (x_{b,i} - x_{a,i}) \end{array} \right] \leq \displaystyle \frac{\partial f}{\partial x_i} \Big |_{x (\tilde \gamma)} (x_{b,i} - x_{a,i} ) \\
\hspace{35mm}\leq \mathop {\max} \left[ \begin{array}{l} \underline \kappa_i (x_{b,i} - x_{a,i}), \\ \overline \kappa_i (x_{b,i} - x_{a,i}) \end{array} \right].
\end{array}
\end{equation}

Substituting this result into (\ref{eq:Lproof4}) then yields (\ref{eq:bound1L}) and (\ref{eq:bound1U}).

\section{Derivation of the Quadratic Bounds}

The derivation is similar to that of the linear case, and simply involves taking the Taylor series expansion one degree higher, with

\vspace{-4mm}
\begin{equation}\label{eq:Qproof2}
\begin{array}{l}
\hat f (1) = \hat f (0) + \displaystyle \frac{d \hat f}{d \gamma} \Big |_{0} + \frac{1}{2} \frac{d^2 \hat f}{d\gamma^2} \Big |_{\tilde \gamma} 
\end{array}
\end{equation}

\noindent for some $\tilde \gamma \in (0,1)$. Applying the chain rule

\vspace{-4mm}
\begin{equation}\label{eq:Qproof3}
\displaystyle \frac{d \hat f}{d \gamma} \Big |_{\gamma} = \displaystyle \mathop {\sum} \limits_{i = 1}^{n} \frac{\partial f}{\partial x_i} \Big |_{x (\gamma)} \frac{d x_i}{d \gamma} \Big |_{\gamma} = \nabla f(x(\gamma))^T (x_b - x_a)
\end{equation}

\noindent and then differentiating once more with respect to $\gamma$ yields

\begin{equation}\label{eq:Qproof4}
\begin{array}{l}
\displaystyle \frac{d^2 \hat f}{d \gamma^2} \Big |_{\gamma} = \displaystyle \mathop {\sum} \limits_{i = 1}^{n} \frac{d}{d \gamma} \left( \frac{\partial f}{\partial x_i} \Big |_{x (\gamma)} \frac{d x_i}{d \gamma} \Big |_{\gamma} \right) \\
\hspace{20mm}= \displaystyle \mathop {\sum} \limits_{i = 1}^{n} \frac{d}{d \gamma} \left( \frac{\partial f}{\partial x_i} \Big |_{x (\gamma)} \right) \frac{d x_i}{d \gamma} \Big |_{\gamma}\;,
\end{array}
\end{equation}

\noindent where we have ignored the terms corresponding to $d^2 x_i / d \gamma^2$ as all such terms are 0. Applying the chain rule again yields

\vspace{-4mm}
\begin{equation}\label{eq:Qproof5}
\begin{array}{l}
\displaystyle \frac{d^2 \hat f}{d \gamma^2} \Big |_{\gamma} = \displaystyle \mathop {\sum} \limits_{i = 1}^{n} \mathop {\sum} \limits_{j = 1}^{n}  \frac{\partial^2 f}{\partial x_j \partial x_i} \Big |_{x (\gamma)} \frac{d x_j}{d \gamma} \Big |_{\gamma} \frac{d x_i}{d \gamma} \Big |_{\gamma} \\
\hspace{10mm} = \displaystyle \mathop {\sum} \limits_{i = 1}^{n} \mathop {\sum} \limits_{j = 1}^{n}  \frac{\partial^2 f}{\partial x_j \partial x_i} \Big |_{x (\gamma)} (x_{b,i} - x_{a,i})(x_{b,j} - x_{a,j}).
\end{array}
\end{equation}

Substituting the results of (\ref{eq:Qproof3}) and (\ref{eq:Qproof5}) into (\ref{eq:Qproof2}), noting that $\hat f (0) = f(x_a)$ and $\hat f (1) = f(x_b)$, and rearranging then leads to 

\vspace{-4mm}
\begin{equation}\label{eq:Qproof6}
\begin{array}{l}
f (x_b) - f (x_a) = \nabla f(x_a)^T (x_b-x_a) + \\
\hspace{10mm}\displaystyle \frac{1}{2} \mathop {\sum} \limits_{i = 1}^{n} \mathop {\sum} \limits_{j = 1}^{n}  \frac{\partial^2 f}{\partial x_j \partial x_i} \Big |_{x (\tilde \gamma)} (x_{b,i} - x_{a,i})(x_{b,j} - x_{a,j}).
\end{array}
\end{equation}

The bounds on the quadratic term are derived in a manner analogous to what was done in the linear case:

\vspace{-4mm}
\begin{equation}\label{eq:mvalbound2}
\begin{array}{l}
(x_{b,i} - x_{a,i})(x_{b,j} - x_{a,j}) \geq 0 \Leftrightarrow \vspace{2mm} \\
\hspace{5mm} \underline M_{ij} (x_{b,i} - x_{a,i})(x_{b,j} - x_{a,j}) \leq \vspace{2mm} \\
\hspace{10mm}\displaystyle \frac{\partial^2 f}{\partial x_j \partial x_i} \Big |_{x(\tilde \gamma)} (x_{b,i} - x_{a,i} )(x_{b,j} - x_{a,j}) \leq \vspace{2mm}\\
\hspace{35mm}\overline M_{ij} (x_{b,i} - x_{a,i})(x_{b,j} - x_{a,j}), \vspace{2mm}\\
(x_{b,i} - x_{a,i})(x_{b,j} - x_{a,j}) \leq 0 \Leftrightarrow \vspace{2mm}\\
\hspace{5mm} \overline M_{ij} (x_{b,i} - x_{a,i})(x_{b,j} - x_{a,j}) \leq \vspace{2mm}\\
\hspace{10mm}\displaystyle \frac{\partial^2 f}{\partial x_j \partial x_i} \Big |_{x(\tilde \gamma)} (x_{b,i} - x_{a,i} )(x_{b,j} - x_{a,j}) \leq \vspace{2mm}\\
\hspace{35mm}\underline M_{ij} (x_{b,i} - x_{a,i})(x_{b,j} - x_{a,j}),
\end{array}
\end{equation}

\noindent and, taking both cases into account, we obtain

\vspace{-4mm}
\begin{equation}\label{eq:mvalbound2pw}
\begin{array}{l}
\mathop {\min} \left[ \begin{array}{l} \underline M_{ij} (x_{b,i} - x_{a,i})(x_{b,j} - x_{a,j}), \\ \overline M_{ij} (x_{b,i} - x_{a,i})(x_{b,j} - x_{a,j}) \end{array} \right] \leq \vspace{2mm} \\
\hspace{10mm}\displaystyle \frac{\partial^2 f}{\partial x_j \partial x_i} \Big |_{x(\tilde \gamma)} (x_{b,i} - x_{a,i} )(x_{b,j} - x_{a,j}) \leq \vspace{2mm}\\
\hspace{20mm}\mathop {\max} \left[ \begin{array}{l} \underline M_{ij} (x_{b,i} - x_{a,i})(x_{b,j} - x_{a,j}), \\ \overline M_{ij} (x_{b,i} - x_{a,i})(x_{b,j} - x_{a,j}) \end{array} \right],
\end{array}
\end{equation}

\noindent which may be substituted into (\ref{eq:Qproof6}) to yield (\ref{eq:bound2L}) and (\ref{eq:bound2U}).

\section{Other Versions}

The bounds (\ref{eq:bound1L})-(\ref{eq:bound2U}) allow a good degree of flexibility by considering both lower and upper bounds on the different partial derivatives. Such flexibility may be useful in certain engineering contexts, where \emph{a priori} knowledge about the system in consideration may be used coherently with the lower and upper bounds on the derivatives \cite{Bunin:13}. However, there are also contexts where these bounds may be needed for purely conceptual reasons and where simpler versions are desired. For example, one might want to suppose \cite{Bunin:13b}:

\vspace{-4mm}
\begin{equation}\label{eq:symm}
\begin{array}{l}
\kappa_i = \overline \kappa_i = - \underline \kappa_i \;, \vspace{2mm} \\
M_{ij} = \overline M_{ij} = -\underline M_{ij} \;,
\end{array}
\end{equation}

\noindent which, if we follow the same steps as before, yields

\vspace{-4mm}
\begin{equation}\label{eq:bound1Ls}
f(x_b) - f(x_a) \geq - \displaystyle \sum_{i=1}^n  \kappa_i | x_{b,i} - x_{a,i} |,
\end{equation}

\vspace{-4mm}
\begin{equation}\label{eq:bound1Us}
f(x_b) - f(x_a) \leq \displaystyle \sum_{i=1}^n  \kappa_i | x_{b,i} - x_{a,i} |,
\end{equation}

\vspace{-4mm}
\begin{equation}\label{eq:bound2Ls}
\begin{array}{l}
f(x_b) - f(x_a) \geq \nabla f (x_a)^T (x_b - x_a) - \\
\hspace{10mm}\displaystyle \frac{1}{2} \sum_{i=1}^n \sum_{j=1}^n  M_{ij} | (x_{b,i} - x_{a,i})(x_{b,j} - x_{a,j}) |,
\end{array}
\end{equation}

\vspace{-4mm}
\begin{equation}\label{eq:bound2Us}
\begin{array}{l}
f(x_b) - f(x_a) \leq \nabla f (x_a)^T (x_b - x_a) + \\
\hspace{10mm}\displaystyle \frac{1}{2} \sum_{i=1}^n \sum_{j=1}^n  M_{ij} | (x_{b,i} - x_{a,i})(x_{b,j} - x_{a,j}) |.
\end{array}
\end{equation}

One may take this one step further and define the bounds with respect to some standard norms. Defining

\vspace{-4mm}
\begin{equation}\label{eq:maxkap}
\kappa = \mathop {\max} \limits_{i=1,...,n} \kappa_i,
\end{equation}

\noindent the bounds (\ref{eq:bound1Ls}) and (\ref{eq:bound1Us}) become

\vspace{-4mm}
\begin{equation}\label{eq:bound1Ls2}
f(x_b) - f(x_a) \geq - \kappa \| x_{b} - x_{a} \|_1,
\end{equation}

\vspace{-4mm}
\begin{equation}\label{eq:bound1Us2}
f(x_b) - f(x_a) \leq \displaystyle \kappa \| x_{b} - x_{a} \|_1.
\end{equation}

For Bounds (\ref{eq:bound2Ls}) and (\ref{eq:bound2Us}), we may consider the following derivation:

\vspace{-4mm}
\begin{equation}\label{eq:Mup}
\begin{array}{l}
\displaystyle \sum_{i=1}^n \sum_{j=1}^n  M_{ij} | (x_{b,i} - x_{a,i})(x_{b,j} - x_{a,j}) | \\
\displaystyle \leq \sum_{i=1}^n \sum_{j=1}^n  M_{ij} | x_{b,i} - x_{a,i}||x_{b,j} - x_{a,j} | \\
\displaystyle \leq \sum_{i=1}^n \sum_{j=1}^n  M_{ij} ( x_{b,i} - x_{a,i})^2,
\end{array}
\end{equation}

\noindent which, with

\vspace{-4mm}
\begin{equation}\label{eq:maxM}
M = \mathop {\max} \limits_{i=1,...,n} \sum_{j=1}^n  M_{ij},
\end{equation}

\noindent allows for (\ref{eq:bound2Ls}) and (\ref{eq:bound2Us}) to be simplified to:

\vspace{-4mm}
\begin{equation}\label{eq:bound2Ls2}
f(x_b) - f(x_a) \geq \nabla f (x_a)^T (x_b - x_a) - \frac{1}{2} M \| x_{b} - x_{a} \|_2^2,
\end{equation}

\vspace{-4mm}
\begin{equation}\label{eq:bound2Us2}
f(x_b) - f(x_a) \leq \nabla f (x_a)^T (x_b - x_a) + \frac{1}{2} M \| x_{b} - x_{a} \|_2^2.
\end{equation}

It may also be shown that the bounds (\ref{eq:bound1L}), (\ref{eq:bound1U}), (\ref{eq:bound2L}), (\ref{eq:bound2U}), (\ref{eq:bound1Ls}), (\ref{eq:bound1Us}), (\ref{eq:bound2Ls}), (\ref{eq:bound2Us}), (\ref{eq:bound1Ls2}), (\ref{eq:bound1Us2}), (\ref{eq:bound2Ls2}), and (\ref{eq:bound2Us2}) all hold with \emph{strict} inequality whenever $x_a \neq x_b$. This follows from (\ref{eq:mvalbound1}) and (\ref{eq:mvalbound2}).

We also refer the reader to \cite{Cartis:13} for more alternatives.

\section{Local Bounds}

As derived, the presented bounds are valid for any arbitrary pair $x_a ,x_b \in \mathcal{X}$, which follows from the validity of the Lipschitz constants over all of $\mathcal{X}$. In certain applications, this globality may, however, add unnecessary conservatism and thus motivate local relaxations \cite{Bunin:13}. Noting that the derivations of the bounds only require them to be valid on the line between $x_a$ and $x_b$, let us define the local Lipschitz constants with respect to these two points in particular as

\vspace{-4mm}
\begin{equation}\label{eq:lip1loc}
\underline \kappa_i^{a,b} < \frac{\partial f}{\partial x_i} \Big |_{x} < \overline \kappa_i^{a,b}, \hspace{10mm} i = 1,...,n, \;\;\; \forall x \in \mathcal{X}_{a,b},
\end{equation}

\vspace{-4mm}
\begin{equation}\label{eq:lip2loc}
\underline M_{ij}^{a,b} < \frac{\partial^2 f}{\partial x_j \partial x_i} \Big |_{x} < \overline M_{ij}^{a,b}, \hspace{3mm} i,j = 1,...,n, \;\;\; \forall x \in \mathcal{X}_{a,b},
\end{equation}

\noindent with

\vspace{-4mm}
\begin{equation}\label{eq:locspace}
\mathcal{X}_{a,b} = \{  x_a + \gamma (x_b - x_a) : \gamma \in [0,1]  \}.
\end{equation}

This then yields the corresponding local versions of (\ref{eq:bound1L})-(\ref{eq:bound2U}):

\vspace{-4mm}
\begin{equation}\label{eq:bound1Lloc}
f(x_b) - f(x_a) \geq \displaystyle \sum_{i=1}^n \mathop {\min} \left[ \begin{array}{l} \underline \kappa_i^{a,b} (x_{b,i} - x_{a,i}), \\ \overline \kappa_i^{a,b} (x_{b,i} - x_{a,i}) \end{array} \right],
\end{equation}

\vspace{-4mm}
\begin{equation}\label{eq:bound1Uloc}
f(x_b) - f(x_a) \leq \displaystyle \sum_{i=1}^n \mathop {\max} \left[ \begin{array}{l} \underline \kappa_i^{a,b} (x_{b,i} - x_{a,i}), \\ \overline \kappa_i^{a,b} (x_{b,i} - x_{a,i}) \end{array} \right],
\end{equation}

\vspace{-4mm}
\begin{equation}\label{eq:bound2Lloc}
\begin{array}{l}
f(x_b) - f(x_a) \geq \nabla f (x_a)^T (x_b - x_a) + \\
\hspace{10mm}\displaystyle \frac{1}{2} \sum_{i=1}^n \sum_{j=1}^n \mathop {\min} \left[ \begin{array}{l} \underline M_{ij}^{a,b} (x_{b,i} - x_{a,i})(x_{b,j} - x_{a,j}), \\ \overline M_{ij}^{a,b} (x_{b,i} - x_{a,i})(x_{b,j} - x_{a,j}) \end{array} \right],
\end{array}
\end{equation}

\vspace{-4mm}
\begin{equation}\label{eq:bound2Uloc}
\begin{array}{l}
f(x_b) - f(x_a) \leq \nabla f (x_a)^T (x_b - x_a) + \\
\hspace{10mm}\displaystyle \frac{1}{2} \sum_{i=1}^n \sum_{j=1}^n \mathop {\max} \left[ \begin{array}{l} \underline M_{ij}^{a,b} (x_{b,i} - x_{a,i})(x_{b,j} - x_{a,j}), \\ \overline M_{ij}^{a,b} (x_{b,i} - x_{a,i})(x_{b,j} - x_{a,j}) \end{array} \right].
\end{array}
\end{equation}



\end{document}